\documentclass[a4paper,abstracton]{scrartcl}
\usepackage[utf8]{inputenc}
\usepackage{amsfonts,amsthm,amssymb,amsmath}
\usepackage[inline]{enumitem}
\usepackage{color}
\usepackage{graphicx} 
\usepackage{authblk}

\renewcommand{\leq}{\leqslant}
\renewcommand{\geq}{\geqslant}
\newtheorem{thm}{Theorem}

\newtheorem{prob}[thm]{Problem}

\theoremstyle{remark}

\DeclareMathOperator{\col}{Col}

\title{Factorizations of regular graphs of infinite degree}
\author{Marcin Stawiski \\ stawiski@agh.edu.pl}

\affil{AGH University of Science and Technology,\\ Faculty of Applied Mathematics, \protect\\al. Mickiewicza 30, 30-059 Krakow, Poland}

\begin{document}

\maketitle
\begin{abstract}
Let $\mathcal{H}=(H_i\colon   i<\alpha )$ be an indexed family of  graphs for some ordinal number $\alpha$.  A family $\mathcal{G}=(G_i\colon   i<\alpha )$ of edge-disjoint subgraphs of a graph $G$ such that $G_i$ is isomorphic to $H_i$ for every $i<\alpha$ and $\bigcup\{E(G_i)\colon   i<\alpha\}=E(G)$ is called a \emph{$\mathcal{H}$-decomposition} of $G$. A $\mathcal{H}$-decomposition of $G$ such that every element of $\mathcal{H}$ is a spanning subgraph of $G$ is called a \emph{$\mathcal{H}$-factorization} of $G$.
Let $\kappa$ be an infinite cardinal. 
Kőnig showed in 1936  that every $\kappa$-regular graph has a factorization into perfect matchings.
 Andersen and Thomassen   proved in 1980 that every $\kappa$-regular connected graph has a $\kappa$-regular spanning tree. 
We generalize both these results by proving the existence of a factorization of every $\kappa$-regular graph into $\lambda$-regular subgraphs for every non-zero $\lambda\leq \kappa$. We investigate indexed families $\mathcal{T}=(T_i\colon  i<\kappa)$ of graphs without isolated vertices such that every connected $\kappa$-regular graph has a $\mathcal{T}$-factorization.
The full classification of such families is given.


\end{abstract}

\section{Introduction}

The study of matchings and related notions is arguably one of the most popular topics in graph theory.
This includes matchability, factorizations, packings, and decompositions.
 One of the most studied problems in this area is a 1-factorization, which is a decomposition into perfect matchings.
The most natural necessary condition for the existence of a 1-factorization of a given graph  is its regularity. In the case of finite graphs this condition is far from being sufficient, even for very simple classes of graphs. The same applies for infinite locally finite graphs (graphs without vertices of infinite degree). However, it turns out that this condition is indeed sufficient in the case of non-locally-finite graphs, as it was shown by Kőnig \cite{konigbookpage} in 1936.
Similarly, Andersen and Thomassen \cite{Thomassen} proved in 1980  the following theorem. 

\begin{thm}[Andersen, Thomassen \cite{Thomassen}]\label{thm:andersen}
If $\kappa$ is an infinite cardinal, then graph $G$ has a spanning $\kappa$-regular tree if and only if $G$ is $\kappa$-regular. 
\end{thm}

 In this paper we focus on possibly the most general decomposition properties for which the only necessary condition is the regularity of a given connected graph. We prove that regularity is a necessary and sufficient condition for non-locally-finite graphs even for the strongest of these properties. 
This covers the mentioned problems such as matchability, factorizations, packings and decompositions, and generalizes the mentioned results of Kőnig, Andersen and Thomassen. Variants  of well-known conjectures by Ringel \cite{Ringel} from 1963 and Gy\'arf\'as \cite{Gyarfas} from 1977 for non-locally-finite graphs and their further strengthenings are discussed in this paper. The proofs of all these variants for non-locally-finite graphs shall follow from the main theorem of this paper.

Let $\mathcal{H}=(H_i\colon    i<\alpha )$ be an indexed family of  graphs for some ordinal number $\alpha$.
We say that  $\mathcal{H}$ \emph{packs} into graph $G$ if there exists a family $(G_i\colon    i<\alpha )$ of edge-disjoint subgraphs of $G$ such that $G_i$ is isomorphic to $H_i$ for every $i<\alpha$.
If $\mathcal{H}$ packs into $G$ and furthermore $\bigcup\{E(G_i)\colon   i<\alpha\}=E(G)$, then  $\mathcal{G}$ is called a \emph{$\mathcal{H}$-decomposition} of $G$.
  \emph{Factor} of a graph $G$ is a spanning subgraph of $G$.  A $\mathcal{H}$-decomposition of $G$ such that every element of $\mathcal{H}$ is a factor of $G$ is called a \emph{$\mathcal{H}$-factorization} of $G$. If $\lambda$ is a cardinal number, then
 a factorization into $\lambda$-regular subgraphs is simply called a \emph{$\lambda$-factorization}.

 Ringel \cite{Ringel} conjectured that $2n+1$ copies of every finite tree with $n$ edges pack into $K_{2n+1}$.  The most straightforward variant for non-locally finite graphs would be the conjecture stating that for every infinite cardinal number $\kappa$ and every tree $T$ on $\kappa$ vertices there exists a packing of $\kappa$ copies of $T$ into $K_\kappa$. 
 As packing in Ringel's Conjecture is a decomposition, we can demand that packing to be a decomposition of $K_\kappa$. Furthermore, we can consider arbitrary $\kappa$-regular connected graphs instead of $K_\kappa$.
In contrast, Gy\'arf\'as \cite{Gyarfas} conjectured that every family $\mathcal{T}=(T_i\colon    2\leq i\leq n )$ of trees such that $T_i$ has order $i$ packs into $K_{n}$. Again, such packing is a decomposition but, unlike in Ringel's Conjecture, trees in family $\mathcal{T}$ are pairwise non-isomorphic. We can propose a variant of Gy\'arf\'as' Conjecture stating that for every infinite cardinal $\kappa$ every family $\mathcal{T}$ of at most $\kappa$ pairwise non-isomorphic trees of order at most $\kappa$ packs into $K_\kappa$. Again, we can consider arbitrary $\kappa$-regular connected graphs instead of $K_\kappa$ and further demand this packing to be a decomposition.
 
 We can combine all the proposed conjectures and ask for which non-locally-finite connected graph $G$  there exists a $\mathcal{T}$-decomposition for an arbitrary family $\mathcal{T}$ of $\kappa$ trees of order at least two and at most $\kappa$.  A necessary condition is the $\kappa$-regularity of $G$.  In this paper we prove that it is also a sufficient condition.
 The main result of this paper is Theorem \ref{thm:main}, which provides a positive answer to all the mentioned conjectures for non-locally finite graphs. Theorem \ref{thm:main} is even stronger as it includes not only packings and decompositions but also factorizations. Furthermore, it provides a complete classification of indexed families $\mathcal{T}=(T_i\colon  i<\kappa)$ of graphs such that each element of $\mathcal{T}$ has no isolated vertex and every $\kappa$-regular connected graph has a $\mathcal{T}$-factorization.

\begin{thm}\label{thm:main}
Let $\kappa$ be an infinite cardinal and let $\mathcal{T}=(T_i\colon  i<\kappa)$ be an indexed family of graphs without isolated vertices. Every connected $\kappa$-regular graph $G$ has a $\mathcal{T}$-factorization if and only if each element of $\mathcal{T}$ is a forest with $\kappa$ components, each of order at most $\kappa$.
\end{thm}

We can apply Theorem \ref{thm:main} to various problems by setting a suitable family $\mathcal{T}$. To show that every family $\mathcal{T}'$ of $\kappa$ forests of order at most $\kappa$ packs into every $\kappa$-regular graph it is enough to partition $\mathcal{T}'$ into $\kappa$ sets of cardinality $\kappa$ (possibly omitting isolated vertices). Thus, we obtain a family $\mathcal{T}$  of $\kappa$ forests with $\kappa$ components of order at most $\kappa$, for which we can apply Theorem \ref{thm:main}.
The same method may be applied to obtain an arbitrary decomposition into $\kappa$ non-trivial forests of size at most $\kappa$. A $\lambda$-factorization for non-zero $\lambda\leq\kappa$ may be obtained by setting $\mathcal{T}$ to be a family of $\kappa$ forests with $\kappa$ components, each isomorphic to $\lambda$-regular tree.
For $\lambda=1$ we obtain the mentioned result of Kőnig. One can easily deduce the existence of $\kappa$-regular spanning tree from the existence of $\kappa$-regular spanning forest by connecting its components. Therefore, for $\lambda=\kappa$ Theorem \ref{thm:main} gives a strengthening of the Theorem \ref{thm:andersen}.

The problem of graph decomposition is closely related to the colouring number. The \emph{colouring number} $\col(G)$ of a graph $G$ is the least cardinal number $\mu$ for which there exists an enumeration $(v_i\colon 0 < |V(G)|)$ of vertices of $G$ such that each vertex $v_i \in V(G)$ has less than $\mu$ neighbours of smaller indices.
Erdős and Hajnal proved \cite{erdoscolouringnumber} in 1967 that if $\lambda$ is an infinite cardinal, then there exists a decomposition of a graph $G$ into a union of $\lambda$ forests if and only if $\col(G)\leq \lambda^+$. Note that the colouring number of the complete graph on $\kappa$ vertices is equal to $\kappa $. Therefore, if $\kappa $ is an infinite cardinal, then each connected infinite $\kappa $-regular graph has a decomposition into $\lambda$ forests if and only if $\kappa =\lambda$ or $\kappa =\lambda^+$. 
It follows that if $\kappa$ is a limit cardinal, then we cannot replace $\mathcal{T}=(T_i\colon  i<\kappa )$ with an indexed family of smaller cardinality. 
However, if $\lambda^+=\kappa$, then for every connected $\kappa$-regular graph $G$ there exists an indexed family $\mathcal{T}=(T_i\colon  i<\lambda )$ of forests with $\kappa$ components, each of order at most $\kappa$ such that $G$ has a $\mathcal{T}$-factorization. It is unknown to the author if there exists an indexed family $\mathcal{T}=(T_i\colon  i<\lambda )$ of forests such that every $\kappa$-regular graph has a $\mathcal{T}$-factorization.

\begin{prob}\label{prob:prob}
Let $\lambda$ be an infinite cardinal. Is there any indexed family $\mathcal{T}=(T_i\colon  i<\lambda)$ of forests  such that every connected $\lambda^+$-regular graph $G$ has a $\mathcal{T}$-factorization?
\end{prob}

\section{Factorizations and decompositions}

When considering $\mathcal{H}$-decomposition, we always implicitly assume that elements of $\mathcal{H}$ are vertex-disjoint.  We consider each ordinal $\alpha=\{\beta\colon   \beta<\alpha \}$ as a well-ordered set with the standard well-ordering of ordinals. If $\alpha,\beta, \gamma$ are ordinals, then we treat the Cartesian products $\alpha \times \beta$ and $\alpha \times \beta \times \gamma$ as well-ordered sets with lexicographic order induced by the well-ordering of ordinals.  For notions of graph theory and set theory  which are not defined in this paper see \cite{Diestel, Jech}.

Our main goal is to  prove Theorem \ref{thm:main}. We divide its proof into three parts. The first part is the theorem below which shows the necessity of conditions in Theorem \ref{thm:main}. Theorem \ref{thm:neccessary} also applies to the Problem \ref{prob:prob}.

\begin{thm}\label{thm:neccessary}
Let $\kappa$ be an infinite cardinal and let $\mathcal{T}$ be an indexed family of at least two graphs without isolated vertices. If $\kappa$-regular tree has a $\mathcal{T}$-factorization, then each element of $\mathcal{T}$ is a forest with $\kappa$ components, each of order at most $\kappa$.
\end{thm}
\begin{proof}
Graph $G$ does not contain a cycle. Hence, each element of $\mathcal{T}$ does not contain a cycle. Suppose to the contrary that there exists an element $T$ of $\mathcal{T}$ which has less than $\kappa$ components. Let $T'$ be a spanning subgraph of $G$. If $T'\cong T$, then $G-T'$ has the same number of components as $T'$. It follows that there exists a vertex $v \in V(G)$ such that each edge incident to $v$ is contained in $T'$. Therefore, $v$ is an isolated vertex in $G-T'$, and $v$ cannot be a vertex in any factor of $G-T'$ without isolated vertices. Elements of $\mathcal{T}$ do not contain isolated vertices. Therefore, there is no $\mathcal{T}$-factorization of $G$.
\end{proof}
The next part of the proof of Theorem \ref{thm:main} is Theorem \ref{thm:kapparegularfactorization}, which covers the case of  factorizations into $\kappa$ many $\kappa$-regular forests.
\begin{thm}\label{thm:kapparegularfactorization}
Let $\kappa$ be an infinite cardinal.  Connected graph  $G$ has a factorization into $\kappa$ regular forests of degree $\kappa$ if and only if $G$ is $\kappa$-regular.
\end{thm}

\begin{proof}
Assume first that $G$ contains a $\kappa$-regular spanning forest. It follows that $G$ has $\kappa$ vertices, each of degree at least $\kappa$. The size of $G$ is $\kappa$. Hence, each vertex of $G$ has degree at most $\kappa$. This proves the necessity of $\kappa$-regularity of $G$. Therefore, it remains to prove the sufficiency of $\kappa$-regularity of $G$.

Let $F$ be a spanning $\kappa$-regular forest of $G$ which exists by Theorem \ref{thm:andersen}. Let $v_0$ be an arbitrary vertex of $G$. Consider the enumeration $(v_i\colon  i<\kappa)$ of vertices of $G$ such that in the rooted forest $(F,v_0)$ if $v_i$ is a son of $v_j$, then $i>j$. For every $j<\kappa$ we partition the set of sons of $v_j$ in $(F,v_0)$  into sets $X_j^m(t)$ for every $m,t<\kappa$, each of cardinality $\kappa$. Note that every vertex of $G$ except $v_0$ belongs to exactly one set from the family $\{X_j^m(t)\colon j,m,t<\kappa\}$.
In the proof, we construct a family  $\{ C_j^m\colon   j,m<\kappa \}$ satisfying:
\begin{enumerate}[label = \textnormal{(C\arabic*)}]
    \item $C_j^m \subset N(v_j)\cap \{v_i\colon   i>j \}$, \label{itm:C1}
    \item $C_i^m \cap C_j^m=\emptyset$, for $i \neq j$, \label{itm:C2}
        \item $C_j^m\cap C_j^n= \emptyset$, for $m \neq n$, \label{itm:C3}
    \item $|C_j^m| =\kappa$, \label{itm:C4}
    \item  if $v_j v_i \in E(G)$ and $i>j$, then there exists $m<\kappa$ such that $v_i \in C_j^m$. \label{itm:C5} 
\end{enumerate}

First, we describe how to obtain a desired $\kappa$-factorization of $G$ into $\kappa$ many $\kappa$-regular forests  using the family $\{ C_j^m\colon   j,m<\kappa \}$ satisfying the  conditions above. 
We construct a $\kappa$-factorization $\mathcal{F}= (F^m \colon   m < \kappa \}$ by setting $E(F^m)= \{v_j v_i \colon   v_i \in C_j^m, j < i < \kappa     \}$. Index $m$ is related to the $m$-th factor. If $v_j$ is a vertex in a component $F$ of $F^m$, then $C_j^m$ is the set of sons of $v_j$ in $(F,v)$ where $v$ denotes the vertex $v_i \in F$ with the least index $i$.

By the condition \ref{itm:C4}, every graph in $\mathcal{F}$ is a $\kappa$-regular spanning subgraph of $G$.
Let $F^m$ be an element of $\mathcal{F}$, and assume that there exists a cycle in $F^m$. The vertex of the greatest index in this cycle has two neighbours in $F^m$. Hence, it belongs to $C_j^m$ and $C_i^m$ for some distinct $i,j<\kappa$. This contradicts \ref{itm:C2}. It follows that $F^m$ is a forest for every $m<\kappa$.
By the conditions \ref{itm:C3} and \ref{itm:C5} every edge of $G$ appears in exactly one element of $\mathcal{F}$. Therefore, $\mathcal{F}$ is a factorization of graph $G$ into $\kappa$ regular forests of degree $\kappa$.

It remains to construct the family $\{ C_j^m\colon   j,m<\kappa \}$. We construct  sets $\{A_{j}^m\colon  j,m<\kappa\}$ and  $\{B_{j}^m\colon  j,m<\kappa\}$, and then we obtain $C_{j}^m$ by setting $C_j^m=A_{j}^m \cup B_{j}^m$ for every $j,m < \kappa$. We shall construct sets $\{A_{j}^m\colon  j,m<\kappa\}$ and  $\{B_{j}^m\colon  j,m<\kappa\}$ by transfinite induction on $(m,\tau,i) \in \kappa \times \kappa \times \kappa$ with respect to the lexicographic order on $\kappa \times \kappa \times \kappa$.
During step $(m,\tau,i)$ we either assign vertex $v_i$ to $a^m_j(y)$ for some $j,y<\kappa$, we put $v_i$ in $B^m_j$ for some $j<\kappa$, or we proceed to the next step without doing anything. Assigning $v_i$ to $a^m_j(y)$  is equivalent to defining $a^m_j(y)$  as $v_i$. Without loss of generality we can assume that $V(G) \cap \kappa =\emptyset$. Initially, we temporarily assign a different ordinal number less than $\kappa$ to each element of $\{a^m_j(y)\colon   m,j,y<\kappa\}$ but still refer to $a^m_j(y)$ as not defined until some vertex $v_i\in V(G)$ is assigned to it. This is done for technical reasons so that sets containing $a^m_j(y)$ are well-defined throughout the entire proof.

Recall that index $m$ is related to the $m$-th factor.
After executing steps $(m,\tau,i)$ for every $\tau,i<\kappa$ set   $\{a^m_j(y)\colon  y <\kappa\}$ has been defined for every $j<\kappa$, and we can define  $A_{j}^m=\{a^m_j(y)\colon  y <\kappa\}$. During steps $(m,\tau,i)$ for $\tau,i<\kappa$ we define set $B^m_j$ by putting vertices in it. At the start of the induction no vertex lies in $B^m_j$ for every $m,j<\kappa$.

For a fixed triple $(m,\tau,i)$ let $\sigma_\tau^m(i)=0$ if no $a_j^m(y)$ has  been defined, let $\sigma_\tau^m(i)$ be the least ordinal for which there exist $j, y \leq\sigma^m_\tau(i)$ such that $a^m_{j}(y)$ has not been  defined, or let $\sigma_\tau^m(i)=\kappa$ if every element of  $\{a^m_j(y)\colon  j,y <\kappa\}$ has already been defined.
The  parameter above and the family $X_j^m(t)$ shall ensure that every element in the set $\{a^m_j(y)\colon  j,y <\kappa\}$ is defined after executing steps $(m,\tau',i')$ for $\tau',i'<\kappa$. In particular, it shall ensure that $|A^m_j|\leq|C^m_j|=\kappa$. Index $\tau$ in a triple $(m,\tau,i)$ is an auxiliary index which  serves the purpose of considering every vertex $v_i$ multiple times. Throughout the induction every vertex can be assigned to more than one element of $\{a_i^{m'}(y)\colon   i,m',y<\kappa\}$ but it can be assigned to at most one of them for the fixed $m'$.  We consider vertices of $G$ one by one, and we assign  $v_i$ to $a_j^m(y)$ if all the  conditions below are satisfied:

\begin{enumerate}[label = \textnormal{(D\arabic*)}]
\item vertex $v_i$ has not been assigned to  $a_{j'}^m(y')$ nor put to $B_{j'}^m$ for every $j',y'<\kappa$, \label{itm:D1}
\item $v_i$ is a neighbour of $v_j$ and $i>j$,  \label{itm:D2}
\item for every $y'<y$, vertex $a^m_{j}(y')$ has already been defined but $a^m_{j}(y)$ has been not, \label{itm:D3}
\item $v_i\notin A_j^{m'}$ for every $m'<m$, \label{itm:D4}
\item $v_i\notin B_j^{m'}$ for every $m'< m$, \label{itm:D5}
\item $v_i \notin X_j^{m''}(t)$ for every $(m'',t)>(m,\tau)$ and $v_i\notin X^m_{j'}(\tau)$ for every $j'\neq j$,  \label{itm:D6} 
\item $j \leq \sigma_\tau^m(i)$ and $y \leq \sigma_\tau^m(i)$, \label{itm:D7}
\item $j$ is the least index for which conditions \ref{itm:D2}--\ref{itm:D7} are satisfied for some $y<\kappa$. \label{itm:D8}
 \end{enumerate}

If $v_i$ has not been assigned to any $a_j^m(y)$ by the  conditions above (for the fixed $m$), then we consider putting $v_i$ in $B_j^m$ for some $j<\kappa$.  If condition \ref{itm:D1} is satisfied and $j$ is the least index for which conditions \ref{itm:D2}, \ref{itm:D4}, \ref{itm:D5} and \ref{itm:D6} are satisfied, then we put $v_i$ in $B_j^m$. Otherwise, we do nothing and proceed to the next index. Induction on $\tau<\kappa$ simply means that we repeat the procedure above  for every $i \in \kappa$. Similarly, induction on $m<\kappa$ means that we repeat the procedure for every $(\tau, i) \in \kappa \times \kappa$.

After the induction on $(m,\tau,i) \in\kappa \times \kappa \times \kappa$ we define $C_{j}^m=A_{j}^m \cup B_{j}^m$ for every $j,m < \kappa$. Note that for every $m,j< \kappa$ sets $A^m_j$ and $B^m_j$ are disjoint. We prove that the family $\{C_j^m\colon  j,m<\kappa\}$, obtained by the recursive construction, satisfy conditions \ref{itm:C1}--\ref{itm:C5}. The first three conditions are easy to check and follow directly from the construction. Condition \ref{itm:C1}, \ref{itm:C2} are satisfied by  conditions \ref{itm:D2} and \ref{itm:D1} respectively in the construction of sets $A_j^m$ and $B_j^m$.  Condition \ref{itm:C3} follows from conditions \ref{itm:D4} and \ref{itm:D5}. We need to show that conditions \ref{itm:C4} and \ref{itm:C5} are satisfied.

Now, we prove that  $|A_j^m| =\kappa$ for every $m,j< \kappa$. This is equivalent to $\{ \sigma^m_\tau(0)\colon   \tau<\kappa \}$ not being bounded by any ordinal less than $\kappa$ for every $m <\kappa$. We show that $\sigma^m_\tau(0)$, as a function of $\tau$, is strictly increasing on the set $\{\tau<\kappa\colon   \sigma^m_\tau(0) <\kappa \}$ for every $m<\kappa$. 

Fix $m< \kappa$. For any non-zero $\alpha<\kappa$ we have $\sup \{\sigma^m_\tau(0)\colon  \tau<\alpha \} \leq \sigma^m_\alpha(0)$. Therefore,  it is enough to prove that for every $\alpha<\kappa$ we have $\sigma^m_\alpha(0) <\sigma^m_{\alpha+1}(0)$ or $\sigma^m_\alpha(0)=\kappa$. Denote $s= \sigma^m_\alpha(0)$.
Before the execution of step $(m,\alpha,0)$ vertex $a_{j}^m(y)$ has already been defined for every $j<s$ and $y<s$. Let $C=\{a^m_s(y)\colon  y\leq s \}\cup \{a^m_j(s)\colon  j\leq s \}$. At least one element of $C$ has not been defined. 
For a fixed $\alpha$  consider the induction on $(m,\alpha,i)$. Notice that for every  $a_j^m(y) \in C$ we have $\kappa$ indices $i$ such that vertex $v_i \in V(G)$  satisfies conditions \ref{itm:D1}, \ref{itm:D2}, \ref{itm:D4}, \ref{itm:D5} and \ref{itm:D6} when it is considered at step $i$.
Indeed, each  element of $X^m_j(\alpha)$ satisfies these conditions for some $a_j^m(y) \in C$. 
If $a_j^m(y) \in C$ has not been assigned before step $i$ and $v_i \in X^m_j(\alpha)$, then $v_i$  satisfy conditions \ref{itm:D1}--\ref{itm:D7} when considered as a candidate for $a_j^m(y')$ for some $y' \leq s$. It follows from the second part of \ref{itm:D6} that condition \ref{itm:D8} is also satisfied. Thus, $v_i=a_j^m(y')$.
Note that there are less than $\kappa$ elements in $C$. Recall that there are $\kappa$ indices $i$ such that vertex $v_i \in V(G)$  satisfies conditions \ref{itm:D1}, \ref{itm:D2}, \ref{itm:D4}, \ref{itm:D5} and \ref{itm:D6} when it is considered at step $i$. Therefore, we assign a vertex for every element of $C$ during the induction on $i$ with fixed $m$ and $\alpha$. 
Thus, $\sigma^m_\alpha(0) <\sigma^m_{\alpha+1}(0)$. We proved that for every $j,m<\kappa$ we have $|A_j^m| =\kappa$. As $A_j^m\subseteq C_j^m$, we obtained that the family $\{C_j^m\colon  j,m<\kappa\}$ satisfies \ref{itm:C4}.

It remains to prove that condition \ref{itm:C5} holds for  $\{C_j^m\colon  j,m<\kappa\}$. Assume that $\{C_j^m\colon  j,m<\kappa\}$ does not satisfy condition \ref{itm:C5}.
Let $(i,j)$ be the least element in $\kappa \times \kappa$  such that $i>j$ and $v_{j}v_{i}\in E$ but $v_i \notin C_j^m$ for every $m<\kappa$. Notice that there exists an index $m''<\kappa$ such that $v_i \notin X_{j'}^m(\tau)$ for every $ m''\leq m$, $\tau < \kappa$, $j'<\kappa$. It follows from the  paragraph above that for every $m < \kappa$ there exists an index $t(m)$ such that $\sigma^m_{t(m)}(i)\geq \max \{i,j \}=i$. For $m''<m$ and $\tau\geq t(m)$ we consider step $(m,\tau,i)$, and we check which of the conditions \ref{itm:D1}-\ref{itm:D8} would be satisfied for assigning $v_i$ to $a^m_j(y)$ for some $y\leq \sigma^m_{t(m)}(i)$. 
 
 Conditions \ref{itm:D4}--\ref{itm:D7} and condition \ref{itm:D2} are satisfied for such choice of $j$ and $y$. However, it may happen that condition \ref{itm:D1}, \ref{itm:D3} or \ref{itm:D8} fails. If condition \ref{itm:D1} fails, then it also fails for every successive step $(m,t',i)$ within $m$. Furthermore, condition \ref{itm:D3} may be satisfied for at most one $y$ and by the assumption we did not put $v_i$ to $B_j^m$.

 When assigning $v_i$ to $a^m_j(y)$, condition \ref{itm:D2} is satisfied only for $j<i$, hence for at most $|i|<\kappa$ indices $j$. Therefore, the satisfaction of condition \ref{itm:D1} depends on only these indices $j'$ for which $j'<i$. Therefore, for all but at most $|i|$ indices $m$ condition \ref{itm:D1} is satisfied at $(m,t,i)$ for every $t < \kappa$. Take $m$ such that $m > m''$ and condition \ref{itm:D1} is satisfied at $(m,t,i)$ for every $t < \kappa$. It means that $v_i$ is not assigned to any element of $A^m_{j'}$ or put in  $B^m_{j'}$ for every $j'<\kappa$ during the induction on $(m,t,i)$ for the fixed $m$ and $i$. Consider the assigning of $v_i$ to $a^m_j(y)$. By the assumption conditions \ref{itm:D4} and \ref{itm:D5} are satisfied. By the choice of $m$ conditions \ref{itm:D1} and \ref{itm:D6} are satisfied. Furthermore, $j$ is the only (and therefore the least) index for which all the conditions \ref{itm:D2}, \ref{itm:D4}, \ref{itm:D5}, \ref{itm:D6} are satisfied.
 Therefore, $v_i$ in step $(m,t,i)$ is assigned to an element of $A^m_j$ or $v_i$ is put in $B^m_j$, which contradicts the assumption.
\end{proof}

The next theorem allows us to further factorize $\kappa$-regular forests from Theorem \ref{thm:kapparegularfactorization}. For an arbitrary graph $H$ denote by $S_H(v,d)$ the \emph{sphere} of radius $d$ and centre $v$ in graph $H$. Similarly, denote by $B_H(v,d)$ the \emph{ball} of radius $d$ and centre $v$ in graph $H$.

\begin{thm}\label{thm:factorization}
Let $\kappa$ be an infinite cardinal. If $\mathcal{T}=(T^m\colon   m<\kappa)$ is an indexed family of forests without isolated vertices and with $\kappa$ components each of order at most $\kappa$, then there exists a  $\mathcal{T}$-factorization of $\kappa$-regular tree.
\end{thm}
\begin{proof}
Denote the $\kappa$-regular tree by $G$. For  $m<\kappa$ let $(t^m_{i}\colon  i<\kappa)$ be an enumeration of vertices of $T^m$.  We shall define  set  $\{y^m_{i}\colon  m,i<\kappa \}$ and graph $Y^m$ for every $m<\kappa$ such that $V(Y^m)=\{y^m_i\colon  i<\kappa\}$, $E(Y^m)=\{y^m_iy^m_j\colon  i,j<\kappa, t^m_i t^m_j\in E(T^m)\}$, and the following conditions shall be satisfied:
\begin{enumerate}[label = \textnormal{(E\arabic*)}]
    \item $f^m\colon  t^m_i\mapsto y^m_i$ is an isomorphism of $T^m$ into $Y^m$ for every $m<\kappa$, \label{itm:E1}
    \item if $xy \in E(G)$, then there exists a unique  $(m,i,j)$ such that $i<j$ and $xy=y^m_i y^m_j$,\label{itm:E2}
    \item $V(Y^m)=V(G)$ for every $m<\kappa$. \label{itm:E3}
\end{enumerate}

Now we show that if conditions \ref{itm:E1}--\ref{itm:E3} hold, then the family $(Y^m\colon   m<\kappa)$ is a $\mathcal{T}$-factorization of $G$. 

Condition \ref{itm:E3} means that each $Y^m$ is a factor of $G$. 
It follows from conditions \ref{itm:E1} and \ref{itm:E2} that $(Y^m\colon   m<\kappa)$ is a $\mathcal{T}$-factorization of $G$. 

Pick any $v_0 \in V(G)$ as a root of $G$. For every $m<\kappa$ we define $y^m_0=v_0$. First, for every $m<\kappa$ we partition the family of components of $T^m$ into sets $T^m_d$ for every $d<\omega$ in such a way that $T^m_0$ is a singleton of component containing $t^m_0$ and $|T^m_d|=\kappa$ for every non-zero $d<\omega$. Furthermore, for a component $T$ of $T^m$ denote by $x_T$ the vertex  $t^m_i$ with the least index $i$ in $T$.
For induction on $d \in \omega$, assume that we already assigned elements in $B_G(v_0,d)$ to some elements of $\{y^m_i\colon  i<\alpha\}$ in such a way that the following conditions are satisfied:
\begin{enumerate}[label = \textnormal{(F\arabic*)}]
    \item every vertex in $B_G(v_0,d)$ has been assigned to exactly one vertex in $\{y^m_j\colon  j<\kappa\}$ for every $m<\kappa$, and only vertices in $B_G(v_0,d)$ have been assigned, \label{itm:F1}
    \item if $y^m_i$ and $y^m_j$ have been defined, then $y^m_i y^m_j \in E(G)$ if and only if $t^m_it^m_j \in E(T^m)$, \label{itm:F2}
    \item if $xy$ is an edge in $G$ between two vertices in $B_G(v_0,d)$, then there exists a unique triple $(m,i,j)$ such that $i<j$ and $xy=y^m_i y^m_j$, \label{itm:F3}
    \item  $y^m_i$ has been defined if and only if  $t^m_i \in B_T(x_T,d-d')$ for some $d'\leq d$ and some $ T\in T^m_{d'}$. \label{itm:F4}
\end{enumerate}

For $y\in S_G(v_0,d)$ we define $W_d(y)$ as the set of the vertices $t^m_i\in T^m$ such that  $y=y^m_i$ for some $m,i<\kappa$. For every $y\in S_G(v_0,d)$ we assign each son of $y$ in $G$ to a unique $y^m_j$ such that $t^m_i\in W_d(y)$ is a neighbour of  $t^m_j$ in $T^m$ and $y^m_j$ has not been defined.
Moreover, every possible $y^m_j$ has to be assigned to some son of $y$. Such assignment is possible because $y$ has $\kappa$ sons and if we put $d'=d$ in condition \ref{itm:F4}, we obtain that there are $\kappa$ possible $y^m_j$ which we can assign to each son of $y$. 
Let $X^m_{d+1}=\{x_T\colon   T \in T^m_{d+1}\}$. 
For each $t^m_i \in X^m_{d+1}$ we assign vertex $y^m_i$ to some vertex $v$ in $S_G(v_0,d+1)$ such that $v$ has not been defined yet as a $y^m_j$ for every $j<\kappa$.
For a fixed $m$ each vertex $y^m_i$ has to be assigned with a different vertex from $S_G(v_0,d+1)$, and each possible $y^m_i$ has to be assigned.

Now we show that before executing step $d$ conditions \ref{itm:F1}--\ref{itm:F4} are satisfied. Each of these conditions are trivially satisfied before executing step $0$. Assume then that $d>1$. It follows directly from the construction that conditions \ref{itm:F1}, \ref{itm:F2} and \ref{itm:F4} are satisfied. 
Let $y^m_iy^m_j$ be an edge between two vertices in $E(G)$ and assume that $i<j$. Further, assume that $y^m_i \in S_G(v_0,d-1)$ and $y^m_j \in S_G(v_0,d)$. 
Notice that if $y^m_j=y^{m'}_{j'}$ for some $j'<\kappa$, $m' \neq m$, then $t^{m'}_{j'}=x_T$ for some $T \in T^{m'}$. Therefore, no neighbour in $Y^{m'}$ of $y^{m'}_{j'}$ lies in $S_G(v_0,d-1)$. It follows that the condition  \ref{itm:F3} is satisfied.

It remains to prove that $\{y^m_{i}\colon  m,i <\kappa \}$ satisfies conditions \ref{itm:E1}--\ref{itm:E3}. Condition \ref{itm:E1} is satisfied by \ref{itm:F2} and \ref{itm:F4}. Condition \ref{itm:E2} is satisfied by  \ref{itm:F3}. It follows directly from condition \ref{itm:F1} that condition \ref{itm:E3} is satisfied.
\end{proof}

The proof of Theorem \ref{thm:main}  follows easily from Theorems \ref{thm:neccessary}, \ref{thm:kapparegularfactorization}, and \ref{thm:factorization}. Theorem \ref{thm:neccessary} shows the necessity of conditions in Theorem \ref{thm:main}.  By Theorem \ref{thm:kapparegularfactorization} we obtain a factorization $(Y^m\colon  m<\kappa )$ of $G$ into $\kappa$ regular forests of degree $\kappa$. Then we partition $\mathcal{T}$ into an indexed family $(U^m\colon  m<\kappa)$ of sets each of cardinality $\kappa$. For every $m<\kappa$ set $U^m$ is an indexed family of $\kappa$ forests without isolated vertices and with $\kappa$ components, each of order at most $\kappa$. By Theorem \ref{thm:factorization}, there exists a $U^m$-factorization $W^m$ of $Y^m$ for every $m<\kappa$. It follows that $\{W \colon  W \in W^m, m<\kappa\}$ forms a  $\mathcal{T}$-factorization of $G$. \qed

\bibliographystyle{abbrv}
\bibliography{sources.bib}

\begin{thebibliography}{1}

\bibitem{Thomassen}
L.~D. {Andersen} and C.~{Thomassen}.
\newblock {The cover index of infinite graphs}.
\newblock {\em {Aequationes Math.}}, 20:244--251, 1980.

\bibitem{Diestel}
R.~{Diestel}.
\newblock {\em Graph Theory}.
\newblock Springer-Verlag, Berlin Heidelberg, 2017.

\bibitem{erdoscolouringnumber}
P.~{Erdős} and A.~{Hajnal}.
\newblock {On decomposition of graphs}.
\newblock {\em {Acta Math. Acad. Sci. Hung.}}, 18:359--377, 1967.

\bibitem{Gyarfas}
A.~{Gy\'arf\'as} and J.~{Lehel}.
\newblock {\em {Packing trees of different order into $K_n$, \textnormal{in
  Combinatorics. (Proc. Fifth Hungarian Colloq., Keszthely, 1976), Vol. 1, pp.
  463--469, Vol. 18 of Colloq. Math. Soc. J\'anos Bolyai}}}.
\newblock North-Holland, Amsterdam-New York, 1978.

\bibitem{Jech}
T.~{Jech}.
\newblock {\em {Set Theory. Third Millenium Edition, revised and expanded}}.
\newblock {Springer-Verlag}, Berlin Heidelberg New York, 2003.

\bibitem{konigbookpage}
D.~{Kőnig}.
\newblock {\em Theorie der endlichen und unendlichen Graphen}.
\newblock Akademische Verlagsgesellschaft, Leipzig, 1936.

\bibitem{Ringel}
G.~{Ringel}.
\newblock {\em Extremal problems in the theory of graphs. Theory of Graphs and
  its Applications (Proc. Int. Symp. Smolenice 1963)}.
\newblock Czech Acad. Sci., Prague, 1963.

\end{thebibliography}

\end{document}